\documentclass[11pt]{article}
\usepackage{amsmath}
\usepackage{amssymb}
\makeatletter\def\@fnsymbol#1{\@arabic{#1}}\makeatother
\newtheorem{theo}{\hspace*{\parindent}Theorem}
\newtheorem{lemma}{\hspace*{\parindent}Lemma}
\def\O{{\cal O}}

\title{Inequalities and monotonicity of ratios
for generalized hypergeometric function}
\author{D. Karp\footnote{Institute of Applied Mathematics, Vladivostok, Russia,
e-mail:\emph{dmkrp@yandex.ru}}~
 and S.\,M.\,Sitnik\footnote{Voronezh Institute of the Ministry of
Internal Affairs of the Russian Federation,
e-mail:\emph{box2008in@gmail.com}}}
\date{}
\begin{document}
\maketitle

\begin{center}
\parbox{12cm}{
\small\textbf{Abstract.}  We find two-sided inequalities for the
generalized hypergeometric function of the form
${_{q+1}}F_{q}(-x)$ with positive parameters restricted by certain
additional conditions. Both lower and upper bounds agree with the
value of ${_{q+1}}F_{q}(-x)$ at the endpoints of positive
semi-axis and are asymptotically precise at one of the endpoints.
The inequalities are derived from a theorem asserting the monotony
of the quotient of two generalized hypergeometric functions with
shifted parameters.  The proofs hinge on a generalized Stieltjes
representation of the generalized hypergeometric function. This
representation also provides  yet another method to deduce the
second Thomae relation for ${_{3}F_{2}}(1)$ and leads to an
integral representations of ${_{4}F_{3}}(x)$ in terms of the
Appell function $F_3$. In the last section of the paper we list
some open questions and conjectures.}
\end{center}

\bigskip


\bigskip


\bigskip

\paragraph{1. Introduction.} We will use standard notation for generalized hypergeometric function:
\begin{equation}\label{eq:pFqdefined}
{_{q+1}F_q}\left((a_{q+1});(b_q);z\right)={_{q+1}F_q}\left(\left.\!\!\begin{array}{c}(a_{q+1})\\
(b_q)\end{array}\right|z\!\right)
=\sum\limits_{n=0}^{\infty}\frac{(a_1)_n(a_2)_n\cdots(a_{q+1})_n}{(b_1)_n\cdots(b_q)_nn!}z^n,
\end{equation}
where  $(a_q)=a_1,a_2,\ldots,{a_q}$ is the parameter list and
$(a)_n=a(a+1)\cdots(a+n-1)$ is the shifted factorial. The series
converges for $|z|<1$ and can be analytically continued onto the
entire complex plane cut along $[1,\infty]$.  The celebrated
continued fraction of Gauss
\begin{multline}\label{eq:Gauss1}
\frac{{_2F_1}(a,b;c;-x)}{{_2F_1}(a-1,b;c-1;-x)}=\frac{c-1}{c-1}\raisebox{-7pt}{+}\frac{b(c-a)x}{c}\raisebox{-7pt}{+}\frac{a(c-b)x}{c+1}
\raisebox{-7pt}{+}\cdots\raisebox{-7pt}{+}
\\
\raisebox{-7pt}{+}\frac{(b+n)(c-a+n)x}{c+2n}\raisebox{-7pt}{+}\frac{(a+n)(c-b+n)x}{c+2n+1}\raisebox{-7pt}{+}\cdots
\end{multline}
implies on setting $a=1$ the classical representation
\begin{multline}\label{eq:Gauss2}
{_2F_1}(1,b;c;-x)=\frac{c-1}{c-1}\raisebox{-7pt}{+}\frac{b(c-1)x}{c}\raisebox{-7pt}{+}\frac{(c-b)x}{c+1}
\raisebox{-7pt}{+}\cdots\\
\raisebox{-7pt}{+}\frac{(b+n)(c-1+n)x}{c+2n}\raisebox{-7pt}{+}\frac{(n+1)(c-b+n)x}{c+2n+1}\raisebox{-7pt}{+}\cdots
\end{multline}
This continued fraction converges and has positive elements when
$x>0$ and $c>b>1$.  Hence, its even convergents form an increasing
sequence and approximate the value of ${_2F_1}(1,b;c;-x)$ from
below, while the odd convergents form a decreasing sequence and
approximate ${_2F_1}(1,b;c;-x)$ from above (see details in
\cite{Cuyt,Henr,LW}).  Taking the first three terms we get:
\begin{equation}\label{eq:2F1viaGauss}
\frac{1}{1+bx/c}<{_2F_1}(1,b;c;-x)<\frac{c(c+1)+(c-b)x}{c(c+1)+c(b+1)x}<1.
\end{equation}
 Gauss derived his continued fraction
(\ref{eq:Gauss1}) from contiguous relations for the hypergeometric
function. Another way to explain both the continued fraction
(\ref{eq:Gauss2}) and the inequality (\ref{eq:2F1viaGauss}) is
through Euler's integral representation
\begin{equation}\label{eq:Euler}
{_2F_1}(1,b;c;-x)=\frac{\Gamma(c)}{\Gamma(b)\Gamma(c-b)}\int\limits_{0}^{1}\frac{t^{b-1}(1-t)^{c-b-1}}{1+xt}dt,
\end{equation}
which shows that ${_2F_1}(1,b;c;-x)$ is a Stieltjes function, i.e.
a function of the form
\[
f(x)=\int\limits_{0}^{1/R}\frac{d\phi(u)}{1+xu}
\]
with a bounded, nondecreasing function $\phi$ taking infinitely
many values (see \cite[Chapter 5]{BG}).  For such functions both a
continued fraction representation with positive elements and a
two-sided estimates through Pad\'{e} approximants are well known
(see \cite{Baker,BG,Cuyt,Henr}). Indeed, it is easy to check that
\[
\frac{1}{1+bx/c}=1-\frac{b}{c}x+\O(x^2),~~~x\to0,
\]\[
\frac{c(c+1)+(c-b)x}{c(c+1)+c(b+1)x}=1-\frac{b}{c}x+\frac{b(b+1)}{c(c+1)}x^2+\O(x^3),~~~x\to0.
\]
Hence, the fraction on the left-hand side of
(\ref{eq:2F1viaGauss}) is the Pad\'{e} approximant of
${_2F_1}(1,b;c;-x)$ of order $[0/1]$ while  the fraction on the
right-hand side is the Pad\'{e} approximation of
 order $[1/1]$ (see \cite{Baker,BG}).   The sequences of
diagonal and paradiagonal one-point Pad\'{e} approximants form
upper and lower bounds uniformly converging to the Stieltjes
function, i.e.
\begin{equation}\label{eq:2F1Pade}
 [M/M-1]_f<[M-1/M]_f<f<[M/M]_f<[M+1/M-1]_f
\end{equation}
for each positive integer $M$. Moreover, these bounds are best
possible with respect to the given number of power series
coefficients \cite[Theorem~15.2]{Baker}. Inequality
(\ref{eq:2F1viaGauss}) is a particular case when $M=1$. These
results have been recently generalized to multi-point Pad\'{e}
approximants by Gilewicz, Pindor, Telega and Tokarzewski (see
\cite{Gilewicz,TokTel}). These papers deal, however, with diagonal
and superdiagonal approximants. The subdiagonal approximants used
here are not discussed. The convergents of the continued fraction
(\ref{eq:Gauss1}) were considered by Belevitch in \cite{Bel}.
Explicit expressions for two-point Pad\'{e} approximants for
ratios of the Gauss hypergeometric functions, confluent
hypergeometric functions and Bessel functions are found in
\cite{WimpBeck}. Some applications of Pad\'{e} approximants to
inequalities for special functions are discussed in
\cite{Sitnik,SitnikPreprint1}.

The lower bound in (\ref{eq:2F1viaGauss}) is not only
asymptotically precise at $x=0$ but also agrees with
${_2F_1}(1,b;c;-x)$ at $x=\infty$.  One problem with
(\ref{eq:2F1viaGauss}) is that the same, unfortunately, is not
true for the upper bound, which reduces to the constant
$(c-b)/[c(b+1)]$ at $x=\infty$.  If we wish to generalize
(\ref{eq:2F1viaGauss}) to ${_{q+1}F_{q}}$ we are also faced with
the problem that neither an analogue of the Gauss continued
fraction for the general ${_{q+1}F_{q}}$ is known nor a multiple
Euler representation generalizing (\ref{eq:Euler}) has the form of
a Stieltjes function. In this paper we solve both problems and
find two-sided estimates for $f(x)={_{q+1}F_{q}}(1,
(a_q);(b_q);-x)$, where the lower bound is asymptotically precise
at $x=0$, the upper bound is asymptotically precise at $x=\infty$
and both agree with the values of $f(x)$ at the endpoints of
$[0,\infty]$.  Next, we extend our inequalities to
${_{q+1}F_{q}}(\sigma, (a_q);(b_q);-x)$ for some values of
$\sigma$.  These results are derived from a somewhat stronger
statement about the monotony of a special quotient of
hypergeometric functions. Our method is based on a generalized
Stieltjes representation for ${_{q+1}F_{q}}((a_{q+1});(b_q);-x)$
but does not utilize the relationship with Pad\'{e} approximants.
In the last section of the paper we list some open problems and
conjectures.

Inequalities for general ${_{q+1}F_{q}}$ are surprisingly rare in
the literature. Most important results are due to Luke
\cite{Luke1,Luke2} who uses inequalities between diagonal and
sub-diagonal Pad\'{e} approximants for $(1+z)^{-\beta}$ and
repeatedly integrates them with respect to beta-distributions or
Laguerre distributions. See section~2 for some details.  Carlson
studied in \cite{Carlson} some inequalities for $R$-hypergeometric
function, which can be expressed in terms of a Lauricella $F_D$
function and as such is a generalization of ${_2F_1}$ to the
multivariate case. For ${_2F_1}(a,b;c;x)$ his restrictions on the
parameters are: $c>b>0$, $x<1$. See section~3 for detailed
comparison of our results for ${_2F_1}(a,b;c;x)$ with those of
Carlson. Inequalities for ${_{q+1}F_{q}}(1,(a_q);(b_q);x)$, $q>1$
are not considered in Carlson's paper.

Buschman \cite{Buschman} uses determinant representations to obtain two-sided inequalities for the Gauss
hypergeometric function ${_2F_1}(a+n,b;c;x)$  in terms of ${_2F_1}(a,b;c;x)$ for positive parameters and
$x\in(0,1)$.  His results were later improved and extended by Joshi and Arya \cite{JA1,JA2}.

Inequalities of a different nature for $|{_{p}F_{q}}|$ and
$\Re({_{p}F_{q}})$ have been obtained by Jahangiri and Silvia in
\cite{JahSil} for the special case when ${_{p}F_{q}}$ is
subordinate (in the sense of analytic function theory) to the
linear fractional transformation.

The papers \cite{ABRVV,Baricz1,Baricz2,Drag1,PV,SitnikPreprint2}
(and many more found in the references therein) consider
inequalities for the Gauss function ${_2F_1}(a,b;c;x)$, the Kummer
function ${_1F_1}(a;c;x)$, the Bessel function ${_0F_1}(a;x)$ and
their ratios for $x\in(0,1)$ (cf. \cite{Sitnik}). The emphasis in
these papers is on fine properties near the singular point $x=1$
of ${_2F_1}(a,b;c;x)$. Note also that the method for proving the
monotonicity of ratios based on Lemma~2.1 from \cite{PV} cannot be
applied to prove our Theorem~\ref{th:pFq-monoton}. Since our
results are valid also for $x<0$, they complement those from
\cite{ABRVV,Baricz1,Baricz2,PV} for ${_2F_1}$. General ${_pF_q}$
is not discussed in the above papers, except some results in
\cite{SitnikPreprint2} that are applicable when certain upper and
lower parameters of ${_pF_q}$ differ by half-integers.

\paragraph{2. Generalized Stieltjes representation and inequalities for  ${_{q+1}F_q}$.}
We begin with a representation of ${_{q+1}F_q}(z)$ as a
generalized Stieltjes transform.
\begin{lemma}\label{lm:GenStieljes}
For $\Re{b_k}>\Re{a_k}>0$,  $k=1,2,\ldots,q$, and
$|\arg(1+z)|<\pi$ we have
\begin{equation}\label{eq:Frepr}
{_{q+1}F_q}\left(\left.\!\!\begin{array}{c}\sigma,(a_{q})\\
(b_q)\end{array}\right|-z\!\right)
=A_0\!\!\int\limits_{0}^{1}\frac{s^{a_1-1}g((a_q);(b_q);s)ds}{(1+sz)^{\sigma}},
\end{equation}
where
\begin{equation}\label{eq:g-defined}
g((a_q);(b_q);s)=\int\limits_{\Lambda_q(s)}[1-s/(t_2\cdots{t_q})]^{b_1-a_1-1}\prod_{k=2}^{q}t_k^{a_k-a_1-1}(1-t_k)^{b_k-a_k-1}\,dt_2\cdots{dt_q},
\end{equation}
\begin{equation}\label{eq:Lambda}
\Lambda_q(s)=[0,1]^{q-1}\cap\{t_2,\ldots,t_q:~t_2\cdots{t_q}>s\}
\end{equation}
and
\begin{equation}\label{eq:A0}
A_0=\Gamma\left[\!\!\begin{array}{c}(b_q)\\
(a_q),(b_q-a_q)\end{array}\!\!\right]
\equiv\frac{\Gamma(b_1)\cdots\Gamma(b_q)}{\Gamma(a_1)\cdots\Gamma(a_q)\Gamma(b_1-a_1)\cdots\Gamma(b_q-a_q)}.
\end{equation}
The function $g((a_q);(b_q);s)$ is invariant under simultaneous
shifts of all parameters\emph{:}
\begin{equation}\label{eq:g-invariance}
g((a_q)+\delta;(b_q)+\delta;s)=g((a_q);(b_q);s)
\end{equation}
for any complex $\delta$.
\end{lemma}
 \textbf{Proof.} The
multiple Euler integral representation for ${_{q+1}F_q}(z)$ reads
\cite[formula 7.2.3(10)]{Prud3}:
\begin{equation}\label{eq:multEuler}
{_{q+1}F_q}\left(\left.\!\!\begin{array}{c}\sigma,(a_{q})\\
(b_q)\end{array}\right|-z\!\right)=\Gamma\left[\!\!\begin{array}{c}(b_q)\\
(a_q),(b_q-a_q)\end{array}\!\!\right]
\int\limits_{[0,1]^{q}}\frac{\prod_{k=1}^{q}t_k^{a_k-1}(1-t_k)^{b_k-a_k-1}}{(1+t_1t_2\ldots{t_q}z)^{\sigma}}\,dt_1\cdots{dt_q}.
\end{equation}
Integration in (\ref{eq:multEuler}) is over the $q$-dimensional
unit cube, $\Re{b_k}>\Re{a_k}>0$ for all $k=1,2,\ldots,q$, and
$|\arg(1+z)|<\pi$.  Formula (\ref{eq:multEuler}) is obtained by
repeated application of the generalized Euler integral
\cite[formula (2.2.2)]{AAR}
\[
{_{q+1}F_q}\left(\left.\!\!\begin{array}{c}(a_{q}),a_{q+1}\\
(b_{q-1}),b_q\end{array}\right|z\!\right)=
\frac{\Gamma(b_q)}{\Gamma(a_{q+1})\Gamma(b_q-a_{q+1})}\int\limits_{0}^{1}
t^{a_{q+1}-1}(1-t)^{b_q-a_{q+1}-1}{_{q}F_{q-1}}\left(\left.\!\!\begin{array}{c}(a_{q})\\
(b_{q-1})\end{array}\right|zt\!\right)dt,
\]
where the last step with $q=1$ is the standard Euler integral for
${_2F_1}$.

We make the variable change $s=t_1t_2\cdots{t_q}$,
$t_1=s/(t_2\cdots{t_q})$, leaving $t_2,\ldots,t_q$ unchanged. The
Jacobian of such transformation is
\[
J(s,t_2,\ldots,t_q)=\frac{1}{t_2\cdots{t_q}}.
\]
Hence, we get from (\ref{eq:multEuler}):
\begin{multline*}
{_{q+1}F_q}\left(\left.\!\!\begin{array}{c}\sigma,(a_{q})\\
(b_q)\end{array}\right|-z\!\right)=A_0\int\limits_{0}^{1}\frac{ds}{(1+sz)^{\sigma}}
\\
\times\int\limits_{\Lambda_q(s)}[s/(t_2\cdots{t_q})]^{a_1-1}[1-s/(t_2\cdots{t_q})]^{b_1-a_1-1}\prod_{k=2}^{q}t_k^{a_k-1}(1-t_k)^{b_k-a_k-1}\,\frac{dt_2\cdots{dt_q}}{t_2\cdots{t_q}}
\\
=A_0\int\limits_{0}^{1}\frac{s^{a_1-1}ds}{(1+sz)^{\sigma}}\int\limits_{\Lambda_q(s)}[1-s/(t_2\cdots{t_q})]^{b_1-a_1-1}\prod_{k=2}^{q}t_k^{a_k-a_1-1}(1-t_k)^{b_k-a_k-1}\,dt_2\cdots{dt_q},
\end{multline*}
where $\Lambda_q(s)$ and $A_0$ are defined in (\ref{eq:Lambda})
and (\ref{eq:A0}), respectively. Introducing the notation
(\ref{eq:g-defined}) for the inner integral we arrive at formula
(\ref{eq:Frepr}).  The shift invariance (\ref{eq:g-invariance}) is
obvious from the definition of $g((a_q);(b_q);s)$.~~$\square$

The representation (\ref{eq:Frepr}) is a key ingredient in the
proof of the following theorem.
\begin{theo}\label{th:pFq-monoton}
Suppose $\delta>0$, $b_k>a_k>0$, $k=1,\ldots,q$. The function
\begin{equation}\label{eq:f-defined}
f(\sigma,(a_q);(b_q);\delta;x)\equiv
\frac{{_{q+1}F_q}\left(\sigma,(a_{q})+\delta;
(b_q)+\delta;-x\right)}{{_{q+1}F_q}\left(\sigma,(a_{q});(b_q);-x\right)}
\end{equation}
is monotone decreasing if $\sigma>0$ and monotone increasing if
$\sigma<0$ for all $x>-1$ .
\end{theo}

\noindent\textbf{Proof.} Put
\[
A_\delta=\Gamma\left[\!\!\begin{array}{c}(b_q)+\delta\\
(a_q)+\delta,(b_q-a_q)\end{array}\!\!\right].
\]
Then by (\ref{eq:Frepr}) and (\ref{eq:g-invariance}):
\[
f(\sigma,(a_q);(b_q);\delta;x)=\frac{A_\delta\int_{0}^{1}s^{a_1+\delta-1}(1+sx)^{-\sigma}g((a_q)+\delta;(b_q)+\delta;s)ds}
{A_0\int_{0}^{1}s^{a_1-1}(1+sx)^{-\sigma}g((a_q);(b_q);s)ds}
\]\[
=\frac{A_\delta\int_{0}^{1}s^{a_1+\delta-1}(1+sx)^{-\sigma}g((a_q);(b_q);s)ds}
{A_0\int_{0}^{1}s^{a_1-1}(1+sx)^{-\sigma}g((a_q);(b_q);s)ds}.
\]
The statement of the theorem is equivalent to saying that
$f'_x(\sigma,(a_q);(b_q);\delta;x)<0$ when $\sigma>0$ and
$f'_x(\sigma,(a_q);(b_q);\delta;x)>0$ when $\sigma<0$.
Differentiating the definition of $f$ above we see that both
inequalities are equivalent to the single inequality
\begin{equation}\label{eq:Cebyshev}
\int\limits_{0}^{1}q(s)p(s)ds\int\limits_{0}^{1}h(s)p(s)ds<
\int\limits_{0}^{1}q(s)h(s)p(s)ds\int\limits_{0}^{1}p(s)ds,
\end{equation}
where
\[
p(s)=s^{a_1+\delta-1}(1+sx)^{-\sigma}g((a_q);(b_q);s),
\]\[
q(s)=s^{\delta},~~~~h(s)=\frac{s}{1+xs}.
\]
The function $p(s)$ is positive, while the functions $q(s)$ and
$h(s)$ are monotone increasing for fixed $x>-1$ and $0<s<1$.
Hence, the above inequality is an instance of the Chebyshev
inequality \cite[Chapter IX, formula
(1.1)]{Mitrinovic}.~~$\square$

The value of $f(\sigma,(a_q);(b_q);\delta;x)$ at $x=\infty$ can be
found using the representation of ${_{q+1}F_q}(z)$ in the
neighborhood of the singular point $z=\infty$. Assume for the
moment that no numerator parameters differ by an integer. Then,
according to \cite[formula 7.2.3.77]{Prud3},

\begin{equation}\label{eq:pFqasymp}
{_{q+1}F_q}\left(\left.\!\!\!\begin{array}{c}(a_{q+1})\\
(b_q)\end{array}\right|-z\!\right)\!=
\Gamma\!\!\left[\!\!\begin{array}{c}(b_{q})\\
(a_{q+1})\end{array}\!\!\right]\sum\limits_{k=1}^{q+1}\Gamma\!\left[\!\!\begin{array}{c}a_k, (a_{q+1})'-a_k\\
(b_{q})-a_k\end{array}\!\!\right]\!{_{q+1}F_q}\left(\left.\!\!\!\begin{array}{c}1+a_k-(b_{q}),a_k\\
1+a_k-(a_{q+1})'\end{array}\right|-1/z\!\right)z^{-a_k},
\end{equation}
where the prime at $(a_{q+1})$ means that the term $a_k$ is
excluded from the list.   It follows that for
$\sigma<\min(a_1,a_2,\ldots,a_q)$ we will have
\begin{equation}\label{eq:pFqasymp1}
{_{q+1}F_q}\left(\left.\!\!\begin{array}{c} \sigma, (a_{q})\\
(b_q)\end{array}\right|-z\!\right)=
\Gamma\!\!\left[\!\!\begin{array}{c} (a_{q})-\sigma, (b_{q})\\
(a_{q}),(b_{q})-\sigma\end{array}\!\!\right]z^{-\sigma}+o(z^{-\sigma}),
~~~z\to\infty.
\end{equation}
 Hence,
\begin{equation}\label{eq:fatinfty}
f(\sigma,(a_q);(b_q);\delta;\infty)=\Gamma\!\!\left[\!\!\begin{array}{c}(a_{q}), (a_{q})+\delta-\sigma, (b_{q})-\sigma, (b_{q})+\delta\\
(a_{q})-\sigma, (a_{q})+\delta,(b_{q}),
(b_{q})+\delta-\sigma\end{array}\!\!\right].
\end{equation}
When some of the numerator parameters on the left-hand side of
(\ref{eq:pFqasymp}) differ by an integer, formula
(\ref{eq:pFqasymp}) breaks down and one has to resort to much more
involved \cite[formula 7.2.3.78]{Prud3}. Formula
(\ref{eq:fatinfty}), however, remains valid by continuity
(${_pF_q}/\Gamma[(b_q)]$ is an entire function of its parameters -
see \cite[7.3.2.8]{Prud3}).

Formula (\ref{eq:fatinfty}) and Theorem~\ref{th:pFq-monoton} imply
the following inequalities valid for $x>0$:
\begin{equation}\label{eq:fineq1}
f(\sigma,(a_q);(b_q);\delta;\infty)<f(\sigma,(a_q);(b_q);\delta;x)<1=f(\sigma,(a_q);(b_q);\delta;0),
\end{equation}
 for positive $\sigma$ and
\begin{equation}\label{eq:fineq2}
1<f(\sigma,(a_q);(b_q);\delta;x)<f(\sigma,(a_q);(b_q);\delta;\infty),
\end{equation}
 for negative $\sigma$.

Taking $\delta=1$ in (\ref{eq:fatinfty}) we will have:
\begin{equation}\label{eq:f1atinfty}
f(\sigma,(a_q);(b_q);1;\infty)=\prod\limits_{i=1}^{q}\frac{b_i(a_i-\sigma)}{a_i(b_i-\sigma)}.
\end{equation}
A simple calculation shows that
\[
1-{_{q+1}F_q}(1,(a_q);(b_q);-x)=x{_{q+1}F_q}(1,(a_q)+1;(b_q)+1;-x)\prod\limits_{i=1}^{q}(a_i/b_i).
\]
It follows that
\begin{multline}\label{eq:pFq-ratio}
\frac{1-{_{q+1}F_q}(1,(a_q);(b_q);-x)}{x{_{q+1}F_q}(1,(a_q);(b_q);-x)}
=\frac{{_{q+1}F_q}(1,(a_q)+1;(b_q)+1;-x)}{{_{q+1}F_q}(1,(a_q);(b_q);-x)}\prod\limits_{i=1}^{q}(a_i/b_i)
\\=f(1,(a_q);(b_q);1;x)\prod\limits_{i=1}^{q}(a_i/b_i).
\end{multline}
Combined with (\ref{eq:fineq1}) and
(\ref{eq:f1atinfty}) this gives:
\[
\prod\limits_{i=1}^q\frac{(a_i-1)}{(b_i-1)}<\frac{1-{_{q+1}F_q}(1,(a_q);(b_q);-x)}{x{_{q+1}F_q}(1,(a_q);(b_q);-x)}<\prod\limits_{i=1}^{q}(a_i/b_i).
\]
Formulas (\ref{eq:fineq1}) and (\ref{eq:f1atinfty}) imply that
these bounds are best possible and by Theorem~\ref{th:pFq-monoton}
the function in the middle is monotone.  A simple rearrangement of
the last formula leads to

\begin{theo}\label{th:pFqmain}
For $b_k>a_k>1$, $k=1,\ldots,q$, and $x>0$ the inequality
\begin{equation}\label{eq:pFqmain}
\frac{1}{1+x\prod_{i=1}^{q}(a_i/b_i)}<{_{q+1}F_q}(1,(a_q);(b_q);-x)<\frac{1}{1+x\prod_{i=1}^{q}[(a_i-1)/(b_i-1)]}
\end{equation}
holds true.
\end{theo}
Some comments are here in order.  Since, clearly,
\[
\frac{1}{1+x\prod_{i=1}^{q}(a_i/b_i)}=1-x\prod_{i=1}^{q}\frac{a_i}{b_i}+\O(x^2),~~~x\to{0},
\]
we conclude that the lower bound in (\ref{eq:pFqmain}) is
asymptotically precise at $x=0$ and coincides with the Pad\'{e}
approximant of order $[0/1]$ to ${_{q+1}F_q}(1,(a_q);(b_q);-x)$ at
zero .  Hence, the inequality from below in (\ref{eq:pFqmain}) can
be proved by noticing that ${_{q+1}F_q}(1,(a_q);(b_q);-x)$ is a
Stieltjes function by Lemma~\ref{lm:GenStieljes}.  This allows to
relax the restrictions on the parameters to $b_k>a_k>0$. Moreover,
a result of Gilewicz and Magnus \cite{GilMAgnus} (see also
\cite[Chapter 25, Theorem~2]{Mitrinovic}) implies that
\emph{\textbf{the lower bound in \emph{(\ref{eq:pFqmain})} is
valid for all $x>-1$.}}  We further generalize the lower bound in
the following theorem.
\begin{theo}\label{th:LowerGeneral}
For $b_k>a_k>0$, $k=1,\ldots,q$, $x>-1$ and $\sigma\geq{1}$ the
inequality
\begin{equation}\label{eq:LowerGeneral}
\frac{1}{\left(1+x\prod_{i=1}^{q}(a_i/b_i)\right)^{\sigma}}<{_{q+1}F_q}(\sigma,(a_q);(b_q);-x)
\end{equation}
holds true.
\end{theo}
\textbf{Proof.} The proof is an application of Jensen's inequality
\begin{equation}\label{eq:Jensen}
\varphi\left(\frac{\int\limits_a^{b}f(s)d\mu(s)}{\int\limits_a^{b}d\mu}\right)
\leq\frac{\int\limits_a^{b}\varphi(f(s))d\mu(s)}{\int\limits_a^{b}d\mu}
\end{equation}
valid for convex $\varphi$, integrable $f$ and non-negative
measure $\mu$ \cite[formula (7.15)]{Mitrinovic}. Take
$\varphi(u)=u^\sigma$, $\sigma\geq{1}$, $f(s)=1/(1+sx)$,
$d\mu(s)=s^{a_1-1}g((a_q);(b_q);s)ds$.  By (\ref{eq:Frepr})
\[
\int\limits_0^{1}f(s)d\mu(s)=\frac{1}{A_0}{_{q+1}F_q}\left(1,(a_{q});(b_q);-x\right),
\]
\[
\int\limits_0^{1}d\mu=\int\limits_0^{1}s^{a_1-1}g((a_q);(b_q);s)ds=
\frac{{_{q+1}F_q}\left(1,(a_{q});(b_q);0\right)}{A_0}=\frac{1}{A_0}
\]
and
\[
\int\limits_0^{1}\varphi(f(s))d\mu(s)=\frac{1}{A_0}{_{q+1}F_q}\left(\sigma,(a_{q});(b_q);-x\right).
\]
Hence, (\ref{eq:Jensen}) reads
\[
\left({_{q+1}F_q}\left(1,(a_{q});(b_q);-x\right)\right)^\sigma\leq{_{q+1}F_q}\left(\sigma,(a_{q});(b_q);-x\right).
\]
Combined with the lower bound from (\ref{eq:pFqmain}) this yields
(\ref{eq:LowerGeneral}).  The restrictions on the parameters are
explained in the remark preceding this theorem. ~$\square$

The inequality (\ref{eq:LowerGeneral}) was previously obtained by
Luke in \cite{Luke1} using a completely different method (see
\textbf{Introduction}). His restrictions are $b_k\geq{a_k}>0$,
$k=1,\ldots,q$, $\sigma>0$ and $x>0$ so that
Theorem~\ref{th:LowerGeneral} extends the validity of Luke's
inequality to $-1<x<0$ under the additional restriction
$\sigma>1$. The case $0<\sigma<1$,  $-1<x<0$ remains open. We
conjecture that (\ref{eq:LowerGeneral}) is still true in this
case.  Another curious lower bound due to Luke \cite{Luke1} valid
for $b_k\geq{a_k}>0$, $k=1,\ldots,q$, $x>0$, $0<\sigma\leq{1}$ is
given by
\[
\frac{1}{\left(1+x\sigma\prod_{i=1}^{q}(a_i/b_i)\right)}<{_{q+1}F_q}(\sigma,(a_q);(b_q);-x).
\]

In a similar fashion but under stronger assumptions we can prove a
generalization of the upper bound from (\ref{eq:pFqmain}):
\begin{theo}\label{th:UpperGeneral}
For $b_k>a_k>1$, $k=1,\ldots,q$, $x>0$ and $0<\sigma\leq{1}$ the
inequality
\begin{equation}\label{eq:UpperGeneral}
{_{q+1}F_q}(\sigma,(a_q);(b_q);-x)<\frac{1}{\left(1+x\prod_{i=1}^{q}[(a_i-1)/(b_i-1)]\right)^\sigma}
\end{equation}
holds true.
\end{theo}
\textbf{Proof.} Again apply (\ref{eq:Jensen}) but this time with
$\varphi(u)=u^{1/\sigma}$, $0<\sigma\leq{1}$,
$f(s)=1/(1+sx)^{\sigma}$ and
$d\mu(s)=s^{a_1-1}g((a_q);(b_q);s)ds$. This yields
\[
\left({_{q+1}F_q}(\sigma,(a_q);(b_q);-x)\right)^{1/\sigma}\leq{_{q+1}F_q}(1,(a_q);(b_q);-x).
\]
The combination of this inequality with the upper bound from
(\ref{eq:pFqmain}) results in (\ref{eq:UpperGeneral}).~~$\square$

Formula (\ref{eq:pFqasymp1}) and the relation
\[
\frac{1}{1+x\prod_{i=1}^{q}[(a_i-1)/(b_i-1)]}=\frac{1}{x}\prod_{i=1}^{q}\frac{b_i-1}{a_i-1}+\O(1/x^2),~~~x\to{\infty},
\]
show that  the upper bound in (\ref{eq:pFqmain}) is asymptotically
precise at $x=\infty$.  This, unfortunately, is not true of the
upper bound (\ref{eq:UpperGeneral}).  We conjecture an
asymptotically precise upper bound for $\sigma\geq{1}$ in the last
section of the paper.  We also remark that, unlike lower bounds,
our upper bounds are very different from those of Luke
\cite{Luke1,Luke2}.

\textbf{Remark.} It is interesting to observe that the proofs of
Theorems~\ref{th:LowerGeneral} and \ref{th:UpperGeneral} work for
any generalized Stieltjes function, i.e. any function of the form
\[
f_\sigma(x)=\int\limits_0^{1/R}\frac{d\phi}{(1+xt)^{\sigma}}
\]
with a bounded nondecreasing function $\phi$ taking infinitely
many values and $\sigma>0$.  We have then
\[
[f_1(x)]^\sigma\leq [f_1(0)]^{\sigma-1}
f_\sigma(x),~~\sigma\geq{1},
\]
and
\[
f_{\sigma}(x)\leq
[f_1(0)]^{1-\sigma}[f_1(x)]^\sigma,~~0<\sigma\leq{1}.
\]
Hence, by combining these inequalities with inequality
(\ref{eq:2F1Pade}) we can ''exponentiate'' the known inequalities
between a Stieltjes function and its Pad\'{e} approximants.

\paragraph{3. The case $q=1$.}  In this case we
are able to extend the upper bound from (\ref{eq:UpperGeneral}) to
negative $x$ as follows:
\begin{theo} For $b>a+1>1$, $0<\sigma\leq{1}$ and $-1<x<0$ the inequality
\begin{equation}\label{eq:2F1upperbound}
{_2F_1}(\sigma,a;b;-x)<\frac{1}{\left(1+\frac{a}{b-1}x\right)^\sigma}
\end{equation}
holds true.
\end{theo}
\textbf{Proof}. The function $f$ defined by (\ref{eq:f-defined})
for $q=1$, $\delta=1$ and $\sigma=1$ takes the form
\[
f(1,a;b;1;x)=\frac{{_2F_1}(1,a+1;b+1;-x)}{{_2F_1}(1,a;b;-x)}.
\]
The value $f(1,a;b;1;-1)$ is finite under the assumptions of the
theorem and is found by the Gauss formula for
${_2F_1}(a_1,a_2;b;1)$:
\[
f(1,a;b;1;-1)=\frac{{_2F_1}(1,a+1;b+1;1)}{{_2F_1}(1,a;b;1)}=\frac{b}{b-1}.
\]
By Theorem~\ref{th:pFq-monoton} $f(1,a;b;1;x)$ is monotone
decreasing for $x\in(-1,0)$ and so
\[
f(1,a;b;1;0)=1<f(1,a;b;1;x)<\frac{b}{b-1}=f(1,a;b;1;-1).
\]
According to (\ref{eq:pFq-ratio})
\[
\frac{a}{b}<\frac{1-{_2F_1}(1,a;b;-x)}{x{_2F_1}(1,a;b;-x)}<\frac{a}{b}\left(\frac{b}{b-1}\right)
\]
or
\[
\frac{1}{1+ax/b}<{_2F_1}(1,a;b;-x)<\frac{1}{1+ax/(b-1)}.
\]
Inequality (\ref{eq:2F1upperbound}) is obtained from this estimate
by an application of Jensen's inequality as in the proof of
Theorem~\ref{th:UpperGeneral}.~$\square$

It is interesting to compare these results with those from
\cite{Carlson}. Carlson's  inequality for the case
$b>a\geq{\sigma}>0$ considered here reads
\begin{equation}\label{eq:Carlson}
\max\left\{\frac{(1+x)^{b-a-\sigma}}{(1+x(1-a/b))^{b-\sigma}},(1+ax/b)^{-\sigma}\right\}<{_2F_1}(\sigma,a;b;-x)<(1+x)^{-\sigma{a}/b}.
\end{equation}
Note first that our condition $b>a\geq{\sigma}>0$ is not more
restrictive than $b>\max(a,\sigma)>0$, since one can exchange the
roles of $a$ and $\sigma$. We see that the lower bound  in
(\ref{eq:Carlson}) is an extension (and possibly a refinement due
to the competitive term under $\max$) of our inequality
(\ref{eq:LowerGeneral}) to the values of $\sigma\in(0,1)$ for the
particular case $q=1$.

To compare the upper bounds we note that for $b>a>\sigma>0$
according to (\ref{eq:pFqasymp1})
\[
{_2F_1}(\sigma,a;b;-x)=\frac{1}{x^\sigma}\left[\frac{\Gamma(b)\Gamma(a-\sigma)}{\Gamma(a)\Gamma(b-\sigma)}\right]^{\sigma}+o(x^{-\sigma}),
~~x\to\infty,
\]
and clearly
\[
(1+x)^{-\sigma{a}/b}=\frac{1}{x^{\sigma{a}/b}}(1+\O(1/x)),~~x\to\infty.
\]
Hence, the upper bound in (\ref{eq:Carlson}) is never
asymptotically precise, albeit agrees with the value
$0={_2F_1}(a,b;c;-x)$ at $x=\infty$.  Our bound in
(\ref{eq:pFqmain}) and the conjectured bound in
(\ref{eq:conjecture}) are asymptotically precise and so both are
better than Carlson's bound at least for large $x$.  The upper
bound (\ref{eq:UpperGeneral}) agrees with the main asymptotic term
in order of $x$ but not in the constant, while the upper bound in
(\ref{eq:Carlson})  has the wrong order and so again our bound
(\ref{eq:UpperGeneral}) is better than Carlson's bound at least
for large $x$.

Finally, for negative $x$ the upper bound in
 (\ref{eq:Carlson}) goes to $\infty$ as $x\to{-1}$, while
under our restrictions on parameters ${_2F_1}(\sigma,a;b;-x)$
remains bounded and so  does our bound in (\ref{eq:UpperGeneral}).
Hence our bound is guaranteed to be better around $x=-1$, while
Carlson's bound is more precise when $x$ is close to $0$.

\paragraph{3. The case $q=2$.} In this case we are able to relax the
assumptions on the parameters imposed in
Lemma~\ref{lm:GenStieljes} and Theorems~\ref{th:pFq-monoton},
\ref{th:pFqmain}, \ref{th:LowerGeneral} and \ref{th:UpperGeneral}.
To this end we give an alternative proof of
Lemma~\ref{lm:GenStieljes} which also shows that the function $g$
defined by (\ref{eq:g-defined}) is expressed for $q=2$ in terms of
${_2F_1}$. Representation (\ref{eq:3F2-int2F1}) below is not new,
it is a slightly different form of \cite[formula
2.21.1.26]{Prud3}. However, we include a short proof which
clarifies the source of restrictions on parameters.
\begin{lemma}\label{lm:3F2-int2F1}
Let $\Re(d+e-b-c)>0$,  $\Re{c}>0$, $\Re{b}>0$ and
$|\arg(1-x)|<\pi$, then
\begin{multline}\label{eq:3F2-int2F1}
{_3F_2}(a,b,c;d,e;x)\\
=\frac{\Gamma(d)\Gamma(e)}{\Gamma(b)\Gamma(c)\Gamma(d+e-b-c)}\int\limits_{0}^{1}\frac{t^{b-1}(1-t)^{d+e-b-c-1}}{(1-xt)^a}{_2F_1}(e-c,d-c;d+e-b-c;1-t)dt.
\end{multline}
\end{lemma}

\noindent\textbf{Proof.} Expand $(1-xt)^{-a}$ on the right-hand
side of (\ref{eq:3F2-int2F1}) into a binomial series and integrate
term by term applying the generalized Euler integral \cite[Theorem
2.2.4]{AAR}
\[
\int\limits_{0}^{1}u^{\gamma-1}(1-u)^{\nu-\gamma-1}{_2F_1}(\alpha,\beta;\gamma;ux)du=\frac{\Gamma(\gamma)\Gamma(\nu-\gamma)}{\Gamma(\nu)}{_2F_1}(\alpha,\beta;\nu;x),
\]
where $\Re{\nu}>\Re{\gamma}>0$, $|\arg(1-x)|<\pi$ (for this
formula to be valid at $x=1$ the additional restriction
$\Re(\nu-\alpha-\beta)>0$ must be imposed), and the Gauss formula
\[
{_2F_1}(\alpha,\beta;\nu;1)=\frac{\Gamma(\nu)\Gamma(\nu-\alpha-\beta)}{\Gamma(\nu-\alpha)\Gamma(\nu-\beta)},
\]
valid under the same restriction. This yields:
\begin{multline}
\int\limits_{0}^{1}t^{b+k-1}(1-t)^{d+e-b-c-1}{_2F_1}(e-c,d-c;d+e-b-c;1-t)dt
\\
=\int\limits_{0}^{1}u^{d+e-b-c-1}(1-u)^{b+k-1}{_2F_1}(e-c,d-c;d+e-b-c;u)du
\\
=\frac{\Gamma(d+e-b-c)\Gamma(b+k)}{\Gamma(d+e-c+k)}{_2F_1}(e-c,d-c;d+e-c+k;1)
\\
=\frac{\Gamma(d+e-b-c)\Gamma(b+k)\Gamma(c+k)}{\Gamma(e+k)\Gamma(d+k)}
\end{multline}
which implies (\ref{eq:3F2-int2F1}) on substitution into the
series. ~~$\square$

\noindent\textbf{Remark.} Representation (\ref{eq:3F2-int2F1})
provides another way to derive Thomae's second fundamental
relation for ${_3F_2}(a,b,c;d,e;1)$ \cite[formula 3.2(2)]{Bailey}.
Indeed, apply the connection formula \cite[formula (2.3.13)]{AAR}
\begin{multline*}
{_2F_1}(e-c,d-c;d+e-b-c;1-t)=\frac{\Gamma(d+e-b-c)\Gamma(c-b)}{\Gamma(d-b)\Gamma(e-b)}{_2F_1}(d-c,e-c;b-c+1;t)
\\
+\frac{\Gamma(d+e-b-c)\Gamma(b-c)}{\Gamma(d-c)\Gamma(e-c)}t^{c-b}{_2F_1}(d-b,e-b;c-b+1;t)
\end{multline*}
on the right-hand side of (\ref{eq:3F2-int2F1}), change
$t\to{1-t}$ and apply (\ref{eq:3F2-int2F1}) again to get
\begin{multline*}
{_3F_2}\!\left[\!\!\!\begin{array}{c}a,b,c\\
d,e\end{array}\!\!\!\right]=\frac{\pi\Gamma(d)\Gamma(e)\Gamma(1-a)\Gamma(d+e-b-c)}{\sin(\pi(c-b))\Gamma(b)\Gamma(c)}
\\
\times\left\{\frac{{_3F_2}[1-a,1-c,d+e-a-b-c;d-a-c+1,e-a-c+1]}{\Gamma(d-a-c+1)\Gamma(e-a-c+1)\Gamma(d-b)\Gamma(e-b)}
-\textrm{idem}(b;c)\right\},
\end{multline*}
where the unit argument is omitted for conciseness and
$\textrm{idem}(b;c)$ after an expression means that the preceding
expression is repeated with $b$ and $c$ interchanged.  To obtain
the form given in \cite[formula 3.2(2)]{Bailey}, one needs to
apply Thomae's first fundamental relation
\begin{multline}\label{eq:Thomae1}
{_3F_2}\!\left[\!\!\!\begin{array}{c}d+e-a-b-c,1-a,1-c\\
d-a-c+1,e-a-c+1\end{array}\!\!\!\right]
\\
=\frac{\Gamma(d-a-c+1)\Gamma(e-a-c+1)\Gamma(b)}{\Gamma(d+e-a-b-c)\Gamma(1+b-a)\Gamma(1+b-c)}
{_3F_2}\!\left[\!\!\!\begin{array}{c}b,1+b-e,1+b-d\\1+b-a,1+b-c\end{array}\!\!\!\right].
\end{multline}

Using Lemma~\ref{lm:3F2-int2F1} we can give a version of
Theorem~\ref{th:pFq-monoton} for the ratio of ${_3F_2}$s under
slightly weaker restrictions on parameters.
\begin{theo}\label{th:ratio}
Let $d+e-b-c>0$, $c>0$, $b>0$, $\delta>0$ and
$\min(b,c)<\min(d,e)$. Then the function
\[
f(a,b,c,d,e,\delta;x)\equiv\frac{{_3F_2}(a,b+\delta,c+\delta;d+\delta,e+\delta;-x)}{{_3F_2}(a,b,c;d,e;-x)}
\]
is monotone decreasing if $a>0$ and monotone increasing if $a<0$
for all $x>-1$.
\end{theo}
\textbf{Proof.} Assume without loss of generality $c=\min(b,c)$
(otherwise exchange the roles of $b$ and $c$). This and the
hypotheses of the theorem imply that $d-c>0$ and $e-c>0$. Now
follow the proof of Theorem~\ref{th:pFq-monoton} to get inequality
(\ref{eq:Cebyshev}) with
\[
p(s)=\frac{s^{b-1}(1-s)^{d+e-b-c-1}}{(1+xs)^a}{_2F_1}(e-c,d-c;d+e-b-c;1-s),
\]\[
q(s)=s^{\delta},~~~~h(s)=\frac{s}{1+xs}.
\]
The function $p(s)$ is positive since $e>c$, $d>c$ and
$d+e-b-c>0$, while the functions $q(s)$ and $h(s)$ are monotone
increasing for fixed $x>-1$, $0<s<1$. The result follows by the
Chebyshev inequality as before.~~$\square$

Inequality (\ref{eq:pFqmain}) takes the form
\begin{equation}\label{eq:3F2ineq}
\frac{1}{1+xbc/de}<{_3F_2}(1,b,c;d,e;-x)<\frac{1}{1+x(b-1)(c-1)/(d-1)(e-1)}
\end{equation}
and here is valid under the assumptions of Theorem~\ref{th:ratio}.
A particular case of this inequality has been used in
\cite{KarpSitnik} to obtain an error bound in the asymptotic
expansion of the first incomplete elliptic integral which leads to
very precise two-sided inequalities for this integral.
Theorems~\ref{th:LowerGeneral} and \ref{th:UpperGeneral} are also
valid here under the assumptions of Theorem~\ref{th:ratio}.

\paragraph{4. The case $q=3$.}  In this section we will show that
Lemma~\ref{lm:GenStieljes} leads to new representations of
${_4F_3}$ as a double integral of ${_2F_1}$ or as a single
integral of the Appell function $F_3$. To this end we need to
demonstrate that the function $g((a_3);(b_3);s)$ can be expressed
in terms of the Appell function $F_3$ or as an integral of
${_2F_1}$. Indeed, we have by (\ref{eq:Frepr})
\begin{multline*}
{_4F_3}(\sigma,(a_3);(b_3);-z)
\\
=A_0\!\!\!\int\limits_{0}^{1}\frac{s^{a_1-1}ds}{(1+sz)^{\sigma}}\!\!\!\int\limits_{\Lambda_3(s)}\!\!\![1-s/(t_2t_3)]^{b_1-a_1-1}t_2^{a_2-a_1-1}t_3^{a_3-a_1-1}(1-t_2)^{b_2-a_2-1}(1-t_3)^{b_3-a_3-1}dt_2dt_3,
\end{multline*}
where
\[
A_0=\frac{\Gamma(b_1)\Gamma(b_2)\Gamma(b_3)}{\Gamma(a_1)\Gamma(a_2)\Gamma(a_3)\Gamma(b_1-a_1)\Gamma(b_2-a_2)\Gamma(b_3-a_3)},
\]\[
\Lambda_3(s)=\{t_2,t_3:t_2t_3>s, 0<t_2<1, 0<t_3<1\}.
\]
The double integral (\ref{eq:g-defined}) here can be written as
follows:
\[
g((a_3);(b_3);s)=\int\limits_{s}^{1}t_2^{a_2-b_1}(1-t_2)^{b_2-a_2-1}dt_2\int\limits_{s/t_2}^{1}[t_2t_3-s]^{b_1-a_1-1}t_3^{a_3-b_1}(1-t_3)^{b_3-a_3-1}dt_3.
\]
Make the change of variables
\[
x=\frac{t_3-s/t_2}{1-s/t_2},~~t_3=x(1-s/t_2)+s/t_2,~~~1-t_3=(1-s/t_2)(1-x),~~dt_3=(1-s/t_2)dx,
\]
in the inner integral to get
\begin{multline}\label{eq:g-doublint}
g((a_3);(b_3);s)=\int\limits_{s}^{1}t_2^{a_2-b_1}(1-t_2)^{b_2-a_2-1}dt_2
\\
\times\int\limits_{0}^{1}[t_2x(1-s/t_2)]^{b_1-a_1-1}(x(1-s/t_2)+s/t_2)^{a_3-b_1}((1-s/t_2)(1-x))^{b_3-a_3-1}(1-s/t_2)dx
\\
=s^{a_3-b_1}\!\!\int\limits_{s}^{1}t_2^{b_1+a_2-a_1-a_3-1}(1-t_2)^{b_2-a_2-1}(1-s/t_2)^{b_1+b_3-a_1-a_3-1}dt_2
\!\!\int\limits_{0}^{1}\frac{x^{b_1-a_1-1}(1-x)^{b_3-a_3-1}}{(1+x(t_2/s-1))^{b_1-a_3}}dx
\\
=\frac{\Gamma(b_1-a_1)\Gamma(b_3-a_3)}{\Gamma(b_1+b_3-a_1-a_3)}s^{a_3-b_1}
\\
\times\!\!\int\limits_{s}^{1}t_2^{b_1+a_2-a_1-a_3-1}(1-t_2)^{b_2-a_2-1}(1-s/t_2)^{b_1+b_3-a_1-a_3-1}
{_2F_1}\left(\left.\!\!\begin{array}{c}b_1-a_3,b_1-a_1\\
b_1+b_3-a_1-a_3\end{array}\right|1-\frac{t_2}{s}\!\right)dt_2
\\
=\frac{\Gamma(b_1-a_1)\Gamma(b_3-a_3)}{\Gamma(b_1+b_3-a_1-a_3)}s^{a_3-a_1}
\\
\times\int\limits_{s}^{1}t_2^{a_2-a_3-1}(1-t_2)^{b_2-a_2-1}(1-s/t_2)^{b_1+b_3-a_1-a_3-1}
{_2F_1}\left(\left.\!\!\begin{array}{c}b_3-a_1,b_1-a_1\\
b_1+b_3-a_1-a_3\end{array}\right|1-\frac{s}{t_2}\!\right)dt_2,
\end{multline}
where we used Euler's integral
\[
\int\limits_{0}^{1}\frac{x^{b_1-a_1-1}(1-x)^{b_3-a_3-1}}{(1+x(t_2/s-1))^{b_1-a_3}}dx=\frac{\Gamma(b_1-a_1)\Gamma(b_3-a_3)}{\Gamma(b_1+b_3-a_1-a_3)}
{_2F_1}\left(\left.\!\!\begin{array}{c}b_1-a_3,b_1-a_1\\b_1+b_3-a_1-a_3\end{array}\right|1-\frac{t_2}{s}\!\right)
\]
and Pfaff's transformation
\[
{_2F_1}\left(\left.\!\!\begin{array}{c}b_1-a_3,b_1-a_1\\b_1+b_3-a_1-a_3\end{array}\right|1-\frac{t_2}{s}\!\right)
=(t_2/s)^{a_1-b1}{_2F_1}\left(\left.\!\!\begin{array}{c}b_3-a_1,b_1-a_1\\b_1+b_3-a_1-a_3\end{array}\right|1-\frac{s}{t_2}\!\right).
\]
Finally, make the change of variables
\[
y=\frac{1-s/t_2}{1-s},~~t_2=\frac{s}{1-y(1-s)},~~~1-t_2=\frac{(1-s)(1-y)}{1-y(1-s)},~~dt_2=\frac{s(1-s)dy}{(1-y(1-s))^2}
\]
in (\ref{eq:g-doublint}) to get:
\begin{multline}\label{eq:g4F3}
g((a_3);(b_3);s)=\frac{\Gamma(b_1-a_1)\Gamma(b_3-a_3)}{\Gamma(b_1+b_3-a_1-a_3)}s^{a_2-a_1}(1-s)^{b_1+b_2+b_3-a_1-a_2-a_3-1}
\\
\times\int\limits_{0}^{1}\frac{y^{b_1+b_3-a_1-a_3-1}(1-y)^{b_2-a_2-1}}{(1-y(1-s))^{b_2-a_3}}
{_2F_1}\left(\left.\!\!\begin{array}{c}b_3-a_1,b_1-a_1\\
b_1+b_3-a_1-a_3\end{array}\right|y(1-s)\!\right)dy.
\end{multline}
This can be further expressed in terms of Appell's function $F_3$,
defined by \cite[formula 5.7(8)]{HTF1}
\[
F_3(\alpha_1,\alpha_2;\beta_1,\beta_2;\gamma;z_1,z_2)
=\sum\limits_{k,l=0}^{\infty}\frac{(\alpha_1)_k(\alpha_2)_l(\beta_1)_k(\beta_2)_l}{(\gamma)_{k+l}}\frac{z_1^k}{k!}\frac{z_2^l}{l!}.
\]
Using the double integral representation \cite[formula
5.8(3)]{HTF1}
\begin{multline}
F_3(\alpha_1,\alpha_2;\beta_1,\beta_2;\gamma;z_1,z_2)=\frac{\Gamma(\gamma)}{\Gamma(\beta_1)\Gamma(\beta_2)\Gamma(\gamma-\beta_1-\beta_2)}
\\
\times\int\limits_{0}^{1}\frac{u^{\beta_1-1}}{(1-uz_1)^{\alpha_1}}\int\limits_{0}^{1-u}\frac{v^{\beta_2-1}(1-u-v)^{\gamma-\beta_1-\beta_2-1}}{(1-vz_2)^{\alpha_2}}dv
\end{multline}
we can obtain by the substitution $t=v/(1-u)$ in the inner
integral, an application of the Euler integral for ${_2F_1}$ and
changing $u\to{y=1-u}$:
\[
\int\limits_{0}^{1}\frac{y^{\gamma-\beta_1-1}(1-y)^{\beta_1-1}}{(1-yz_1)^{\alpha_1}}{_2F_1}\left(\left.\!\!\begin{array}{c}\alpha_2,\beta_2\\
\gamma-\beta_1\end{array}\right|yz_2\!\right)dy=\frac{\Gamma(\gamma)(1-z_1)^{-\alpha_1}}{\Gamma(\beta_1)\Gamma(\gamma-\beta_1)}
F_3(\alpha_1,\alpha_2;\beta_1,\beta_2;\gamma;z_1/(z_1-1),z_2).
\]
This formula is a slightly different guise of \cite[formula
2.21.1.20]{Prud3}.  Combined with (\ref{eq:g4F3}) it yields
\begin{multline}\label{eq:gF3}
g((a_3);(b_3);s)=
\frac{\Gamma(b_1-a_1)\Gamma(b_3-a_3)\Gamma(b_1+b_2+b_3-a_1-a_2-a_3)}{\Gamma(b_2-a_2)[\Gamma(b_1+b_3-a_1-a_3)]^2}
s^{a_2+a_3-a_1-b_2}
\\
\times(1-s)^{b_1+b_2+b_3-a_1-a_2-a_3-1}F_3(b_2-a_3,b_3-a_1;b_2-a_2,b_1-a_1;b_1+b_2+b_3-a_1-a_2-a_3;1-1/s,1-s).
\end{multline}
Finally, we obtain the following representations for ${_4F_3}$:
\begin{multline}\label{eq:4F3-iint2F1}
{_4F_3}(\sigma,(a_3);(b_3);-z)=\frac{\Gamma(b_1)\Gamma(b_2)\Gamma(b_3)}{\Gamma(a_1)\Gamma(a_2)\Gamma(a_3)\Gamma(b_2-a_2)\Gamma(b_1+b_3-a_1-a_3)}
\\
\times\!\!\int\limits_{0}^{1}\!\!\!\int\limits_{0}^{1}\frac{s^{a_2-1}y^{b_1+b_3-a_1-a_3-1}(1-s)^{b_1+b_2+b_3-a_1-a_2-a_3-1}(1-y)^{b_2-a_2-1}}
{(1+sz)^{\sigma}(1-y(1-s))^{b_2-a_3}}
{_2F_1}\!\!\left[\left.\!\!\begin{array}{c}b_3-a_1,b_1-a_1\\
b_1+b_3-a_1-a_3\end{array}\!\right|y(1-s)\!\right]\!\!dsdy,
\end{multline}
\begin{multline}\label{eq:4F3-intF3}
{_4F_3}(\sigma,(a_3);(b_3);-z)=B_1\int\limits_{0}^{1}\frac{s^{a_2+a_3-b_2-1}(1-s)^{b_1+b_2+b_3-a_1-a_2-a_3-1}}{(1+sz)^{\sigma}}
\\
\times
F_3(b_2-a_3,b_3-a_1;b_2-a_2,b_1-a_1;b_1+b_2+b_3-a_1-a_2-a_3;1-1/s,1-s)ds,
\end{multline}
where
\[
B_1=\frac{\Gamma(b_1)\Gamma(b_2)\Gamma(b_3)\Gamma(b_1+b_2+b_3-a_1-a_2-a_3)}{\Gamma(a_1)\Gamma(a_2)\Gamma(a_3)[\Gamma(b_2-a_2)]^2[\Gamma(b_1+b_3-a_1-a_3)]^2}.
\]
Both representations (\ref{eq:4F3-iint2F1}) and
(\ref{eq:4F3-intF3}) are believed to be new.  They have been
verified by termwise integration of the series expansions of
hypergeometric functions occurring in the integrands and comparing
coefficients at powers of $(-z)$.  According to \cite[formula
7.2.4.74)]{Prud3} the function $F_3(1-1/s,1-s)$ encountered in
(\ref{eq:4F3-intF3}) can be expressed as the sum of three
${_3F_2}(s)$.

\paragraph{5. Open questions and conjectures.}
The cases $q=1,2,3$ leave little doubt  that the function
$g((a_q);(b_q);s)$ defined by (\ref{eq:g-defined}) can be
expressed in terms of multiple hypergeometric functions for all
$q$.  The restrictions on parameters in Lemma~\ref{lm:GenStieljes}
are rooted in the definition of $g$ as the integral
(\ref{eq:g-defined}). These observations motivate

\textbf{Open problem~1:} How to express the function
$g((a_q);(b_q);s)$ defined by (\ref{eq:g-defined}) in terms of
multivariate hypergeometric functions for $q>3$? How to extend the
validity of Lemma~\ref{lm:GenStieljes} to a wider range of
parameter values using the analytic continuation of
$g((a_q);(b_q);s)$?

Numerical experiments show that the condition $b_k>a_k>0$ is
sufficient but in no way necessary to the validity of
Theorem~\ref{th:pFq-monoton}.  Hence, our

\textbf{Open problem~2:} How to relax the restrictions on
parameters in Theorem~~\ref{th:pFq-monoton}?

Numerical tests also indicate clearly the following

\textbf{Conjecture~1.} Theorem~\ref{th:LowerGeneral} is true for
all $\sigma>0$ and $\sum_{i=1}^{q}(b_i-a_i)>0$.

The asymptotic formula (\ref{eq:pFqasymp1}) and numerical
experiments suggest the following

\textbf{Conjecture~2.} For $0<\sigma\leq\min(a_1,a_2,\ldots,a_q)$,
$\sum_{i=1}^{q}(b_i-a_i)>0$ and $x>0$:
\begin{equation}\label{eq:conjecture}
{_{q+1}F_q}(\sigma,(a_q);(b_q);-x)
<\frac{1}{\left(1+x\prod_{i=1}^{q}\frac{\Gamma(a_{i})\Gamma(b_{i}-\sigma)}{\Gamma(b_{i})\Gamma(a_{i}-\sigma)}\right)^\sigma}.
\end{equation}

Combined with (\ref{eq:3F2-int2F1}) Thomae's first relation
(\ref{eq:Thomae1}) leads to the following curious identity
\[
\int\limits_{0}^{1}t^{b-1}(1-t)^{d+e-a-b-c-1}{_2F_1}(e-c,d-c;d+e-b-c;1-t)dt
\]\[
=
\frac{\Gamma(b)\Gamma(c)\Gamma(d+e-a-b-c)}{\Gamma(a)\Gamma(d-a)\Gamma(e-a)}
\int\limits_{0}^{1}t^{e-a-1}(1-t)^{a-1}{_2F_1}(e-c,e-b;d+e-b-c;1-t)dt.
\]

\textbf{Open problem~3:} How can one derive the above identity
directly from the properties of ${_2F_1}$? Such derivation would
immediately give another proof for the first fundamental relation
of Thomae.

\paragraph{6. Acknowledgements.} The first author is supported by
the Russian Basic Research Fund (grant no.\,08-01-00028-a), Far
Eastern Branch of the Russian Academy of Sciences  (grant no.\,
06-III-B-01-020) and Presidential grant for leading scientific
schools (grant 2810.2008.1)


\begin{thebibliography}{99}
\bibitem{ABRVV}G.D.\,Anderson, R.W.\,Barnard, K.C.\,Richards,
M.K.\,Vamanamurthy and M.\,Vuorinen, Inequalities for
zero--balanced hypergeometric functions, \emph{Trans.
Amer.Math.Soc.}, 347 (1995), 1713-1723.
\bibitem{AAR} G.E.\,Andrews, R.\,Askey and R.\,Roy, \emph{Special
functions}, Cambridge University Press, 1999.
\bibitem{Bailey} W.N.\,Bailey, \emph{Generalized hypergeometric
series}, Cambridge University Press, 1935.
\bibitem{Baker} G.A.\,Baker,  \emph{Essentials of Pad\'{e} Approximants}, Academic Press, London,
1975.
\bibitem{BG} G.A.\,Baker Jr. and P.\,Graves-Morris, \emph{Pad\'{e}
approximants. Part~I: Basic Theory}, Addison-Wesley Publishing
Company, 1981.
\bibitem{Baricz1}A.\,Baricz, Functional inequalities involving special
functions, \emph{J. Math. Anal. Appl.}, \textbf{319}, 2(2006),
450-459.
\bibitem{Baricz2}A.\,Baricz, Functional inequalities involving special
functions II,  \emph{J. Math. Anal. Appl.}, \textbf{327}, 2(2007),
1202-1213.
\bibitem{Bel}V.\,Belevitch, The Gauss hypergeometric ratio as a positive real
function, \emph{SIAM J. Math. Anal.}, \textbf{13}, 6 (1982),
1024-1040.
\bibitem{Buschman} R.G.\,Buschman, Inequalities for
hypergeometric functions, \emph{Mathematics of computation},
\textbf{30}, 134 (1976), 303-305.
\bibitem{Carlson} B.C.\,Carlson, Some inequalities for
hypergeometric functions, \emph{Proc. of Amer. Math. Soc.}, vol. 17,
no.1 (1966), 32-39.
\bibitem{Cuyt}A.\,Cuyt, V.B.\,Petersen, B.\,Verdonk, H.\,Waadeland, W.B.\,Jones, \emph{Handbook of Continued Fractions for Special Functions}, Springer, 2008.
\bibitem{Drag1}S.S.\,Dragomir, P.\,Cerone, \emph{Advances in Inequalities for Special Functions}, Nova Science, 2007.
\bibitem{Henr} P.\,Henrici, \emph{Applied and computational complex analysis, Vol.2}, Wiley, 1977.
\bibitem{HTF1} A.\,Erd\'{e}lyi, W.\,Magnus, F.~Oberhettinger and F.G.~Tricomi,
{\em Higher transcendental functions}, Vol. 1, McGraw-Hill Book
Company, Inc., New York, 1953.
\bibitem{GilMAgnus}J.\,Gilewicz and A.P.\,Magnus, Sharp
inequalities for the Pad\'{e} approximant errors in the Stieltjes
case, \emph{Rocky Math. J.}, \textbf{21} (1991), 227-233.
\bibitem{Gilewicz}J.\,Gilewicz, M.\,Pindor, J.J.\,Telega and S.\,Tokarzewski,
N-point Pad\'{e} approximants and two-sided estimates of errors on
the real axis for Stieltjes functions, \emph{J. of Comp. and Appl.
Math.} \textbf{178} (2005), 247-253.
\bibitem{JahSil} M.\,Jahangiri and E.M.\,Silvia, Some
inequalities involving generalized hypergeometric functions,
In:\emph{Univalent functions, fractional calculus, and their
applications}, H.M.\,Srivastava and S.\,Owa (eds.), John Wiley and
Sons, 1989.
\bibitem{JA1}  C.M.\,Joshi and J.P.\,Arya, Inequalities for Certain Hypergeometric
Functions,  \emph{Mathematics of Computation}, Vol. 38, No. 157(1982), 201-205.
\bibitem{JA2}  C.M.\,Joshi and J.P.\,Arya, Some inequalities for the Gauss and Kummer hypergeometric
functions, \emph{Indian J. Pure Appl. Math.} 22, No.8 (1991),
637-644.
\bibitem{KarpSitnik} D.\,Karp and S.M.\,Sitnik,
Asymptotic approximations for the first incomplete elliptic
integral near logarithmic singularity, \emph{J. of Comp. and Appl.
Math.}, 205 , No.1 (2007), 186-206,
http://dx.doi.org/10.1016/j.cam.2006.04.053
\bibitem{LW}L.\,Lorentzen and H.\,Waadeland, \emph{Continued fractions with
applications}, North-Holland, 1992.
\bibitem{Luke1} Y.\,L.\,Luke, Inequalities for generalized
hypergeometric functions, \emph{Journal of Approximation Theory},
\textbf{5} (1972), 41-65.
\bibitem{Luke2} Y.\,L.\,Luke, \emph{The Special Functions and Their
Approximations}, Volume II, Academic Press, 1969.
\bibitem{Mitrinovic}D.S.\,Mitrinov\'{c}, J.E.\,Pecari\'{c}, A.M.\,Fink, \emph{Classical and new
inequalities in Analysis}, Kluwer Academic Publishers, 1993.
\bibitem{PV}S.Ponnusamy and M.Vuorinen, Asymptotic expansions and
inequalities for hypergeometric functions, \emph{Mathematika},
\textbf{44} (1997), 278-301.
\bibitem{Prud3}A.P.\,Prudnikov, Yu.A.\,Brychkov, O.I.\,Marichev, \emph{Integrals and
series, vol.3, Additional Chapters}, Moscow Nauka, 1986. English
translation: \emph{Integrals and series. Volume 3: More special
functions.} New York: Gordon and Breach Science Publishers, 1990.
\bibitem{Sitnik} S.M.\,Sitnik,  Inequalities for Bessel functions, \emph{Dokl. Math.} 51 (1995), No.1, 25-28;
translation from \emph{Dokl. Akad. Nauk, Ross. Akad. Nauk} 340 (1995), No.1, 29-32.
\bibitem{SitnikPreprint1}L.A.\,Minin, S.M.\,Sitnik, \emph{Pade approximants for elementary and special functions},
preprint, Institute of Automation and Control Processes, Rus. Acad. Sci., Vladivostok, 1991.
\bibitem{SitnikPreprint2}S.M.\,Sitnik, \emph{Inequalities for the complete Legendre elliptic integrals},
preprint, Institute of Automation and Control Processes, Rus.
Acad. Sci., Vladivostok, 1994.
\bibitem{TokTel} S.\,Tokarzewski and J.J.\,Telega, Inequalities for Two-Point Pad\'{e} Approximants
to the Expansions of Stieltjes Functions in a Real Domain,
\emph{Acta Applicandae Mathematicae}, \textbf{48} (1997), 285-297.
\bibitem{WimpBeck}J.\,Wimp and B.\,Beckermann,
Families of two-point Pad\'{e} approximants and some ${}_4 F_3
(1)$ identities, \emph{SIAM J. Math. Anal.} \textbf{26}, 3 (1995),
761-773.
\end{thebibliography}
\end{document}